\documentclass[12pt]{amsart}%{article}

\usepackage{amsmath,amssymb,amsthm}

\begin{document}

\newcommand{\non}{\nonumber}
\newcommand{\scl}{\scriptstyle}
\newcommand{\scloplus}{{\scl\bigoplus}}
\newcommand{\ot}{\otimes}
\newcommand{\tss}{\hspace{1pt}}
\newcommand{\ts}{\,}
\newcommand{\Y}{{\rm Y}}
\newcommand{\CC}{\mathbb{C}}
\newcommand{\ZZ}{\mathbb{Z}}
\newcommand{\gl}{\mathfrak{gl}}
\newcommand{\al}{\alpha}
\newcommand{\de}{\delta}
\newcommand{\la}{\lambda}
\newcommand{\La}{\Lambda}
\newcommand{\ze}{\zeta}

\newtheorem{thm}{Theorem}[section]
\newtheorem{lemma}[thm]{Lemma}
\newtheorem{cor}[thm]{Corollary}

\newcommand{\bth}{\begin{thm}}
\renewcommand{\eth}{\end{thm}}
\newcommand{\ble}{\begin{lemma}}
\newcommand{\ele}{\end{lemma}}
\newcommand{\bco}{\begin{cor}}
\newcommand{\eco}{\end{cor}}

\newcommand{\bpf}{\begin{proof}}
\newcommand{\epf}{\end{proof}}

\def\beql#1{\begin{equation}\label{#1}}

\newcommand{\bal}{\begin{aligned}}
\newcommand{\eal}{\end{aligned}}
\newcommand{\beq}{\begin{equation}}
\newcommand{\eeq}{\end{equation}}
\newcommand{\ben}{\begin{equation*}}
\newcommand{\een}{\end{equation*}}

\title[Gelfand--Tsetlin bases]
{Gelfand--Tsetlin bases for representations\\[0.2em]
of finite W-algebras and shifted Yangians}

\author[V. Futorny]{Vyacheslav Futorny}
\address{Institute of Mathematics and Statistics\\
University of S\~ao Paulo\\
Caixa Postal 66281- CEP 05315-970\\
S\~ao Paulo, Brazil}
\email{futorny@ime.usp.br}
\author[A. Molev]{Alexander Molev}
\address{
School of Mathematics and Statistics\\
University of Sydney,
NSW 2006, Australia}
\email{alexm@maths.usyd.edu.au}
\author[S. Ovsienko]{Serge Ovsienko}
\address{
Faculty of Mechanics and Mathematics\\
Kiev Taras Shevchenko University\\
Vla\-di\-mir\-skaya 64, 00133, Kiev, Ukraine}
\email{ovsienko@zeos.net}

\date{} % November 2007

\maketitle

\begin{abstract}
Remarkable subalgebras of the Yangian
for $\gl_n$ called the shifted Yangians were
introduced in a recent work by Brundan and Kleshchev
in relation to their study of finite $W$-algebras.
In particular, in that work
a classification of finite-dimensional
irreducible representations of the
shifted Yangians and the associated
finite $W$-algebras was given.
We construct a class of these
representations in an explicit form
via bases of Gelfand--Tsetlin type.
\end{abstract}

\section{Introduction}

A striking relationship between the Yangians and
finite $W$-algebras was first discovered by
Ragoucy and Sorba~\cite{rs:yr}; see also
Briot and Ragoucy~\cite{br:rp}. This relationship
was developed in full generality by
Brundan and Kleshchev~\cite{bk:sy}. The finite $W$-algebras
associated to nilpotent orbits in the Lie algebra $\gl_N$
turned out to be isomorphic to quotients of certain
subalgebras of the Yangian $\Y(\gl_n)$. These subalgebras,
called the {\it shifted Yangians\/} in \cite{bk:sy},
admit a description in terms of generators and relations.
This leads to respective presentations
of the finite $W$-algebras and thus provides new tools to study
their structure and representations. The representation
theory of the shifted Yangians and associated $W$-algebras
was developed in a subsequent paper by
Brundan and Kleshchev~\cite{bk:rs} where
deep connections of the
shifted Yangian representation theory were established.
In particular, a classification of the finite-dimensional
irreducible representations of the shifted Yangians
and the finite $W$-algebras was given in terms
of their highest weights. Moreover,
in the case of the shifted Yangian associated to $\gl_2$
all such representations were explicitly constructed.

Our aim in this paper is to construct in an explicit form
a family of representations of the shifted Yangians
and finite $W$-algebras via bases of Gelfand--Tsetlin type.
Such bases for certain classes of representations
of the Yangian $\Y(\gl_n)$
were constructed in different ways by
Nazarov and Tarasov~\cite{nt:yg, nt:ry}
and Molev~\cite{m:gt}.
We mainly employ
the approach of \cite{nt:yg, nt:ry} which
turns out to be more suitable
for the generalization to the case of the shifted Yangians.
In more detail,
following \cite{bk:sy},
consider an $n$-tuple of
positive integers $\pi=(p_1,\dots,p_n)$ such that
$p_1\leqslant\dots\leqslant p_n$. We can visualize
$\pi$ as a {\it pyramid\/} of left-justified rows of bricks,
where the top row contains $p_1$ bricks, the
second row contains $p_2$ bricks, etc.
Such a pyramid determines a finite $W$-algebra
which we denote by $W(\pi)$.
For each $k\in\{1,\dots,n\}$ we let $\pi_k$ denote
the pyramid with the rows $(p_1,\dots,p_k)$.
Our basis is consistent with the chain of subalgebras
\beql{chainw}
W(\pi_1)\subset W(\pi_2) \subset\dots\subset W(\pi_n).
\eeq
In the case of the one-column pyramid $(1,\dots,1)$ of height $n$
we recover the classical Gelfand--Tsetlin basis
for representations of the Lie algebra $\gl_n$.
For any $\pi$,
the formulas for the action
of the Drinfeld generators of $W(\pi)$
in the basis turn out to be
quite similar to the Yangian case.
These explicit constructions
of representations of $W(\pi)$ proved to be useful
for a description of the Harish-Chandra modules
over finite $W$-algebras and a proof of the associated
Gelfand--Kirillov conjecture based on recent
results of Futorny and Ovsienko~\cite{fo:gas, fo:gar}; see our
forthcoming paper \cite{fmo:gk}.

\section{Shifted Yangians and finite $W$-algebras}

As in \cite{bk:sy}, given a pyramid
$\pi=(p_1,\dots,p_n)$ with
$p_1\leqslant\dots\leqslant p_n$, introduce
the corresponding {\it shifted Yangian\/} $\Y_{\pi}(\gl_n)$
as the associative algebra defined by generators
\beql{gener}
\bal
d^{\ts(r)}_i,&\quad i=1,\dots,n,&&\quad r\geqslant 1,\\
f^{(r)}_i,&\quad i=1,\dots,n-1,&&\quad r\geqslant 1,\\
e^{(r)}_i,&\quad i=1,\dots,n-1,&&\quad r\geqslant p_{i+1}-p_i+1,
\eal
\eeq
subject to the following relations:
\begin{align}
[d_{i}^{\ts(r)},d_{j}^{\ts(s)}]&=0,
\non\\
[e_{i}^{(r)},f_{j}^{(s)}]&=-\ts \de_{ij}\ts\sum_{t=0}^{r+s-1}
d_{i}^{\tss\prime\ts(t)}\ts d_{i+1}^{\ts(r+s-t-1)},\non\\
[d_{i}^{\ts(r)},e_{j}^{(s)}]&=(\de_{ij}-\de_{i,j+1})\ts\sum_{t=0}^{r-1}
d_{i}^{\ts(t)}\ts e_{j}^{(r+s-t-1)},\non\\
[d_{i}^{\ts(r)},f_{j}^{(s)}]&=
(\de_{i,j+1}-\de_{ij})\ts\sum_{t=0}^{r-1}
f_{j}^{(r+s-t-1)}\ts d_{i}^{\ts(t)},
\non
\end{align}
\begin{align}
[e_{i}^{(r)},e_{i}^{(s+1)}]-[e_{i}^{(r+1)},e_{i}^{(s)}]&=
e_{i}^{(r)}e_{i}^{(s)}+e_{i}^{(s)}e_{i}^{(r)},\non\\
[f_{i}^{(r+1)},f_{i}^{(s)}]-[f_{i}^{(r)},f_{i}^{(s+1)}]&=
f_{i}^{(r)} f_{i}^{(s)}+f_{i}^{(s)} f_{i}^{(r)},\non\\
[e_{i}^{(r)},e_{i+1}^{(s+1)}]-[e_{i}^{(r+1)},e_{i+1}^{(s)}]&=
-e_{i}^{(r)}e_{i+1}^{(s)},\non\\
[f_{i}^{(r+1)},f_{i+1}^{(s)}]-[f_{i}^{(r)},f_{i+1}^{(s+1)}]&=
-f_{i+1}^{(s)}f_{i}^{(r)},
\non
\end{align}
\begin{alignat}{2}
[e_{i}^{(r)},e_{j}^{(s)}]&=0\qquad&&\text{if}\quad |i-j|>1,\non\\
[f_{i}^{(r)},f_{j}^{(s)}]&=0\qquad&&\text{if}\quad |i-j|>1,\non\\
[e_{i}^{(r)},[e_{i}^{(s)},e_{j}^{(t)}]]
&+[e_{i}^{(s)},[e_{i}^{(r)},e_{j}^{(t)}]]=0
\qquad&&\text{if}\quad |i-j|=1,\non\\
[f_{i}^{(r)},[f_{i}^{(s)},f_{j}^{(t)}]]
&+[f_{i}^{(s)},[f_{i}^{(r)},f_{j}^{(t)}]]=0
\qquad&&\text{if}\quad |i-j|=1,
\non
\end{alignat}
for all admissible $i,j,r,s,t$, where $d_{i}^{(0)}=1$ and
the elements
$d_{i}^{\tss\prime\ts(r)}$ are found from the relations
\ben
\sum_{t=0}^r d_{i}^{\tss(t)}\ts d_{i}^{\tss\prime\ts(r-t)}=\de_{r0},
\qquad r=0,1,\dots.
\een
Note that the algebra $\Y_{\pi}(\gl_n)$ depends only on
the differences $p_{i+1}-p_i$.
In the particular case
of a rectangular pyramid $\pi$ with $p_1=\dots=p_n$, the
algebra $\Y_{\pi}(\gl_n)$ is isomorphic to the {\it Yangian\/}
$\Y(\gl_n)$; see e.g. \cite{m:yc} for the description
of its structure and representations.
The isomorphism with the $RTT$ presentation of
$\Y(\gl_n)$ was constructed in \cite{bk:pp} providing a proof
of the original result of Drinfeld~\cite{d:nr}.
Moreover, for an arbitrary pyramid $\pi$,
the shifted Yangian $\Y_{\pi}(\gl_n)$
can be regarded as a natural subalgebra of $\Y(\gl_n)$.
Note also that the shifted Yangians can be defined
for more general types of pyramids. However, in accordance to
\cite{bk:sy}, each of these algebras is
isomorphic to $\Y_{\pi}(\gl_n)$
for an appropriate left-justified pyramid $\pi$.

Introduce formal generating series in $u^{-1}$ by
\ben
\bal
d_i(u)&=1+\sum_{r=1}^{\infty} d_i^{\ts(r)}\ts u^{-r},\qquad
f_i(u)=\sum_{r=1}^{\infty} f_i^{(r)}\ts u^{-r},\\
e_i(u)&=\sum_{r=p_{i+1}-p_i+1}^{\infty} e_i^{(r)}\ts u^{-r}
\eal
\een
and set
\ben
a_i(u)=d_1(u)\ts d_2(u-1)\dots d_i(u-i+1)
\een
for $i=1,\dots,n$,\  and
\ben
b_i(u)=a_i(u)\ts e_i(u-i+1),\qquad
c_i(u)=f_i(u-i+1)\ts a_i(u)
\een
for $i=1,\dots,n-1$. It is clear that the
coefficients of the series $a_i(u)$, $b_i(u)$ and $c_i(u)$
generate the algebra $\Y_{\pi}(\gl_n)$. It is not difficult
to rewrite the defining relations in
terms of these coefficients. We point out
a few of these relations here which will be
frequently used later on;
see also \cite{bk:pp}.
We have
\begin{align}
\label{abc}
[a_i(u), c_j(v)]&=0,\qquad
[b_i(u), c_j(v)]=0,\qquad\text{if}\quad i\ne j,\\
\label{crcr}
[c_i(u), c_j(v)]&=0,\qquad\text{if}\quad |i-j|\ne 1,\\
\label{arcrv}
(u-v)\ts [a_i(u), c_i(v)]&=c_i(u)\ts a_i(v)-
c_i(v)\ts a_i(u).
\end{align}

Let $N$ be the number of bricks in the pyramid $\pi$.
Due to the main result of \cite{bk:sy},
the {\it finite $W$-algebra\/} $W(\pi)$,
associated to $\gl_N$ and the pyramid $\pi$, can be defined
as the quotient of $\Y_{\pi}(\gl_n)$ by the two-sided ideal
generated by all elements $d_{1}^{\tss(r)}$ with $r\geqslant p_1+1$.
We refer the reader to \cite{bk:sy, bk:rs}
for a discussion of the origins of the finite $W$-algebras
and more references. Note that
in the case of a rectangular pyramid of height $p$,
the algebra $W(\pi)$ is isomorphic to
the {\it Yangian of level $p$\/}; this relationship was
originally observed in \cite{br:rp} and \cite{rs:yr}.

We will use the same notation for the images of the elements of
$\Y_{\pi}(\gl_n)$ in the quotient algebra $W(\pi)$. Set
\ben
A_i(u)=u^{p_1}\ts (u-1)^{p_2}\ts\dots (u-i+1)^{p_i}\ts a_i(u)
\een
for $i=1,\dots,n$,\  and
\ben
\bal
B_i(u)&=u^{p_1}\ts (u-1)^{p_2}\ts
\dots (u-i+2)^{p_{i-1}}\ts (u-i+1)^{p_{i+1}}\ts b_i(u),\\
C_i(u)&=u^{p_1}\ts (u-1)^{p_2}\ts\dots (u-i+1)^{p_i}\ts c_i(u)
\eal
\een
for $i=1,\dots,n-1$. The following lemma
is immediate from the results of
Brown and Brundan~\cite{bb:ei}.
Here we regard $A_i(u)$, $B_i(u)$, and $C_i(u)$ as series with
coefficients in $W(\pi)$.

\ble\label{lem:polabc}
All series $A_i(u)$, $B_i(u)$, and $C_i(u)$
are polynomials in $u$.
\ele

\bpf
In terms of the $RTT$ presentation
of the Yangian,
each of the series $a_i(u)\in\Y(\gl_n)[[u^{-1}]]$
coincides with a quantum minor
of the matrix of the generators; see
\cite[Theorem~8.6]{bk:pp}. Therefore
the statement for the $A_i(u)$ follows from
the results of \cite[Section~3]{bb:ei}. Note that
the polynomial $A_i(u)$ in $u$ is monic of degree $p_1+\dots+p_i$.
Furthermore, the defining relations of $\Y_{\pi}(\gl_n)$ imply
$[f_i^{(1)},a_i(u)]=c_i(u)$,
and so $C_i(u)=[f_i^{(1)},A_i(u)]$ is a polynomial in $u$
of degree $p_1+\dots+p_i-1$. Similarly,
\ben
b_i(u)\ts (u-i+1)^{p_{i+1}-p_i}=[a_i(u),e^{(p_{i+1}-p_i+1)}_i],
\een
which gives
\ben
B_i(u)=[A_i(u),e^{(p_{i+1}-p_i+1)}_i],
\een
so that $B_i(u)$ is a polynomial in $u$
of degree $p_1+\dots+p_i-1$.
\epf

Note that by \cite[Theorem~6.10]{bk:rs}, all coefficients
of the polynomial $A_n(u)$
belong to the center of $W(\pi)$ and these coefficients
(excluding the leading one) are algebraically independent
generators of the center.

For $i=1,\dots,n-1$
define the elements $h_i^{\ts(r)}\in \Y_{\pi}(\gl_n)$
by the expansion
\ben
1+\sum_{r=1}^{\infty} h_i^{\ts(r)}\ts u^{-r}
=d_i(u)^{-1}\ts d_{i+1}(u)
\een
and set
\ben
H^{(r)}_{i}(u)=u^r+u^{r-1}\ts
h_i^{\ts(1)}+\dots+h_i^{\ts(r)}.
\een

\ble\label{lem:relbc}
For $i=1,\dots,n-1$
in the algebra $W(\pi)$ we have
\ben
(u-v)\ts [B_i(u),C_i(v)]
=A'_{i+1}(u)\ts A_i(v)-A'_{i+1}(v)\ts A_i(u),
\een
where $A'_{i+1}(u)$ is the polynomial in $u$
with coefficients in $W(\pi)$
given by
\ben
\bal
A'_{i+1}(u)&=u^{p_1}\ts (u-1)^{p_2}\dots \ts
(u-i+2)^{p_{i-1}}\ts(u-i+1)^{p_{i+1}}\\
{}&\times a_{i}(u+1)^{-1}
\big(a_{i+1}(u+1)\ts a_{i-1}(u)+c_{i}(u+1)\ts b_{i}(u)\big)\\
{}&{}-H^{(p_{i+1}-p_i)}_{i}(u-i+1)\ts A_i(u).
\eal
\een
Moreover,
\ben
\bal
B_i(u)\ts C_i(u-1)&=A'_{i+1}(u)\ts A_{i}(u-1)
-A_{i+1}(u)\ts A_{i-1}(u-1)\\
{}&+
H^{(p_{i+1}-p_i)}_{i}(u-i)\ts A_i(u)\ts A_i(u-1).
\eal
\een
\ele

\bpf
Observe that for any fixed $i\in\{1,\dots,n-1\}$ the
elements $d^{\tss(r)}_i$, $d^{\tss(r)}_{i+1}$, $e^{(r)}_i$ and
$f^{(r)}_i$ of $\Y_{\pi}(\gl_n)$ satisfy the defining relations
of the shifted Yangian $\Y_{\pi_i}(\gl_2)$, where
$\pi_i=(p_i,p_{i+1})$.
Therefore,
it suffices to prove the first relation in the case $i=1$;
the proof for the remaining values of $i$ will then
easily follow.
Working in the Yangian $\Y(\gl_2)$, we can derive
the relation
\ben
(u-v-1)\ts [d_1(u),e_1(v)]=\big(e_1(v)-e_1(u)\big)\ts d_1(u);
\een
see e.g. \cite{bk:pp}.
This allows us to calculate the commutators $[d_1(u),e_1^{(r)}]$
and leads to an equivalent expression for $b_1(u)$
in the subalgebra $\Y_{\pi}(\gl_2)$:
\ben
b_1(u)=d_1(u)\ts e_1(u)=(1-u^{-1})^{p_2-p_1}\ts e_1(u-1)\ts d_1(u).
\een
Furthermore, starting from the relations
\ben
[e_{1}^{(r)},f_{1}^{(s)}]=-\sum_{t=0}^{r+s-1}
d_{1}^{\tss\prime\ts(t)}\ts d_{2}^{(r+s-t-1)}
\een
in $\Y_{\pi}(\gl_2)$, it is now straightforward to derive that
\ben
\bal
(u-v)\ts [b_1(u),c_1(v)]=a_{1}(u+1)^{-1}
\big(a_2(u+1)+c_1(u+1)\ts b_1(u)\big)\ts a_1(v)&\\
{}-(u^{-1} v)^{p_2-p_1}
a_{1}(v+1)^{-1}
\big(a_2(v+1)+c_1(v+1)\ts b_1(v)\big)\ts a_1(u)&\\
{}-u^{p_1-p_2}
\big(H^{(p_{2}-p_1)}_{1}(u)-H^{(p_{2}-p_1)}_{1}(v)\big)
\ts a_1(u)\ts a_1(v)&.
\eal
\een
The desired relation in $W(\pi)$ is then obtained by
multiplying both sides
by the product $u^{p_2}\ts v^{p_1}$. Furthermore,
by the defining relations,
\begin{multline}
u^{p_2}\ts a_{1}(u+1)^{-1}
\big(a_2(u+1)+c_1(u+1)\ts b_1(u)\big)\\
{}=u^{p_2}\ts\big(d_2(u)+f_1(u)\ts d_1(u)\ts e_1(u)\big).
\non
\end{multline}
This is a polynomial in $u$ due to \cite[Theorem~3.5]{bk:rs}.
Hence, by Lemma~\ref{lem:polabc},
$A'_2(u)$ is a polynomial in $u$ too.

The second part of the lemma is implied by the first
part by taking into account
the relations in the shifted Yangian $\Y_{\pi}(\gl_n)$,
\ben
a_i(u)^{-1}c_i(u)=c_i(u-1)\ts a_i(u-1)^{-1}
\een
and
\ben
(u-i)^{p_{i+1}-p_i}\ts a_i(u-1)^{-1}\ts b_i(u-1)
=(u-i+1)^{p_{i+1}-p_i}\ts b_i(u)\ts a_i(u)^{-1},
\een
which are implied by the defining relations.
\epf

\section{Construction of basis vectors}

Using the canonical homomorphism $\Y_{\pi}(\gl_n)\to W(\pi)$
we can extend
every representation of the finite $W$-algebra $W(\pi)$
to the shifted Yangian $\Y_{\pi}(\gl_n)$.
In what follows we work with representations of $W(\pi)$,
and the results can be easily interpreted in the
shifted Yangian context.

Let us recall some definitions
and results from \cite{bk:rs} regarding
representations of $W(\pi)$. Fix an $n$-tuple
$\la(u)=\big(\la_1(u),\dots,\la_n(u)\big)$
of monic polynomials in $u$
with coefficients in $\CC$, where $\la_i(u)$
has degree $p_i$. We let $L(\la(u))$ denote
the irreducible highest weight representation
of $W(\pi)$ with the highest weight $\la(u)$.
Then $L(\la(u))$ is generated by
a nonzero vector $\ze$ (the highest vector)
such that
\begin{alignat}{2}
B_{i}(u)\ts\ze&=0 \qquad &&\text{for} \quad
i=1,\dots,n-1, \qquad \text{and}
\non\\
u^{p_i}\ts d_{i}(u)\ts\ze&=\la_i(u)\ts\ze \qquad &&\text{for}
\quad i=1,\dots,n.
\non
\end{alignat}
Write
\ben
\la_i(u)=(u+\la^{(1)}_i)\ts(u+\la^{(2)}_i)
\dots (u+\la^{(p_i)}_i),\qquad i=1,\dots,n.
\een
We will assume that the parameters $\la^{(k)}_i$
satisfy the conditions:
for any value $k\in\{1,\dots,p_i\}$ we have
\ben
\la^{(k)}_i-\la^{(k)}_{i+1}\in\ZZ_+,\qquad i=1,\dots,n-1,
\een
where $\ZZ_+$ denotes the set of nonnegative integers.
In this case the representation $L(\la(u))$ of $W(\pi)$
is finite-dimensional.

Denote by $q_k$ the number
of bricks in the column $k$ of the pyramid $\pi$.
We have
$q_1\geqslant\dots\geqslant q_l>0$, where $l=p_n$
is the number of the columns in $\pi$.
If $p_{i-1}< k\leqslant p_i$ for some $i\in\{1,\dots,n\}$
(taking $p_0=0$), then we set
$\la^{(k)}=(\la^{(k)}_i,\dots,\la^{(k)}_{n})$.
Then $q_k=n-i+1$.
Let $L(\la^{(k)})$ denote
the finite-dimensional irreducible representation
of the Lie algebra $\gl_{q_k}$ with the highest weight
$\la^{(k)}$. The vector space
\beql{tenprp}
L(\la^{(1)})\ot\dots\ot L(\la^{(l)})
\eeq
can be equipped with an action of the algebra $W(\pi)$,
and $L(\la(u))$
is isomorphic to a subquotient of the module
\eqref{tenprp}. In particular,
\beql{comdim}
\dim L(\la(u))\leqslant\prod_{k=1}^l \dim L(\la^{(k)}).
\eeq

In what follows we will only consider a certain
family of representations
of $W(\pi)$ by imposing a
{\it generality condition\/}
on the highest weights of
the representations $L(\la(u))$.
We will assume that
\ben
\la_i^{(k)}-\la_j^{(m)}\notin\ZZ,\qquad\text{for all}\ \
i,j \quad\text{and all}\ \  k\ne m.
\een

The {\it Gelfand--Tsetlin pattern\/}
$\La(u)$ (associated with $\la(u)$) is an array of
monic polynomials in $u$ of the form

\ben
\bal
\la^{}_{n1}(u)\qquad\qquad
\la^{}_{n2}(u)&\qquad\quad \dots\quad\qquad \la^{}_{nn}(u)\\[1em]
\la^{}_{n-1,1}(u)\qquad \dots&\qquad\ \  \la^{}_{n-1,n-1}(u)\\[0.2em]
\dots\qquad&\dots\\[0.2em]
\la^{}_{21}(u)\ \ \quad&\quad\la^{}_{22}(u)\\[1em]
\la^{}_{11}&(u)
\eal
\een

\bigskip
\noindent
where
\ben
\la_{ri}(u)=(u+\la_{ri}^{(1)})\dots(u+\la_{ri}^{(p_i)}),
\qquad 1\leqslant i\leqslant r\leqslant n,
\een
with $\la_{ni}^{(k)}=\la_{i}^{(k)}$
and the following
conditions hold
\ben
\la_{r+1,i}^{(k)}-\la_{ri}^{(k)}\in\ZZ_+\qquad\text{and}\qquad
\la_{ri}^{(k)}-\la_{r+1,i+1}^{(k)}\in\ZZ_+
\een
for $k=1,\dots,p_i$ and
$1\leqslant i\leqslant r\leqslant n-1$. We have
$\la^{}_{ni}(u)=\la_i(u)$ for $i=1,\dots,n$, so that
the top row coincides with $\la(u)$.

Most arguments in the rest of the paper will not be
essentially different from
\cite[Section~3]{nt:ry}, so we only sketch the main steps
in the construction of the basis.
Given a pattern $\La(u)$,
introduce the corresponding element $\ze^{}_{\La}$ of
$L(\la(u))$ by the formula
\begin{multline}
\ze^{}_{\La}=\prod_{i=1,\dots,\ts n-1}^{\longrightarrow}
\Big\{\ts \prod_{k=1}^{p_i}
\Big(\ts C_{n-1}(-l^{\tss(k)}_{n-1,i}-1)\dots
C_{n-1}(-l^{\tss(k)}_i)\Big)\\
{}\times \prod_{k=1}^{p_i}
\Big(\ts C_{n-2}(-l^{\tss(k)}_{n-2,i}-1)
\dots C_{n-2}(-l^{\tss(k)}_i+1)\ts C_{n-2}(-l^{\tss(k)}_i)\Big)\\
{}\times \dots \times \prod_{k=1}^{p_i}
\Big(\ts C_{i}(-l^{\tss(k)}_{ii}-1)
\dots C_{i}(-l^{\tss(k)}_i+1)\ts C_{i}(-l^{\tss(k)}_i)\Big)\Big\}
\ts\ze,
\non
\end{multline}
where we have used the notation
\ben
l^{\tss(k)}_i=\la^{(k)}_i-i+1\qquad\text{and}\qquad
l^{\tss(k)}_{ri}=\la^{(k)}_{ri}-i+1.
\een
Note that by \eqref{crcr} we have $[C_i(u),C_i(v)]=0$,
so that the order of the factors in the products over $k$
is irrelevant.

\ble\label{lem:arze}
We have
\ben
A_r(u)\ts \ze^{}_{\La}
=\la_{r1}(u)\dots\la_{rr}(u-r+1)\ts \ze^{}_{\La},
\een
for $r=1,\dots,n$.
\ele

\bpf
When applying $A_r(u)$ to $\ze^{}_{\La}$, separating the first
factor, we need to calculate
$A_r(u)\ts C_s(v)\ts\eta$ for the respective value of $v$.
By \eqref{abc}, the operator
$A_r(u)$ commutes with $C_s(v)$ for $s\ne r$.
Furthermore, by \eqref{arcrv},
\ben
A_r(u)\ts C_r(v)\ts\eta=
\frac{1}{u-v}\ts C_r(u)\ts A_r(v)\ts\eta
+\frac{u-v-1}{u-v}\ts C_r(v)\ts A_r(u)\ts\eta.
\een
The calculation is completed by
induction on the number
of factors $C_i(v)$ in the expression for $\ze^{}_{\La}$,
taking into account that $A_r(v)\ts\eta=0$.
\epf

\ble\label{lem:brze}
For any $1\leqslant i\leqslant r\leqslant n-1$
and $k=1,\dots,p_i$ we have
\ben
\bal
B_r(-l^{\tss(k)}_{ri})\ts \ze_{\La}=
{}&-\la_{1}(-l^{\tss(k)}_{ri})
\dots\la_{i}(-l^{\tss(k)}_{ri}-i+1)\\
{}&\times\la_{r+1,i+1}(-l^{\tss(k)}_{ri}-i)
\dots\la_{r+1,r+1}(-l^{\tss(k)}_{ri}-r)\\
{}&\times\la_{1}(-l^{\tss(k)}_{ri}-1)
\dots\la_{i-1}(-l^{\tss(k)}_{ri}-i+1)\\
{}&\times\la_{r-1,i}(-l^{\tss(k)}_{ri}-i)
\dots\la_{r-1,r-1}(-l^{\tss(k)}_{ri}-r+1)
\ts \ze_{\La+\de_{ri}^{\tss(k)}},
\eal
\een
where $\ze^{}_{\La+\de_{ri}^{(k)}}$
corresponds to the pattern
obtained from $\La(u)$ by replacing $\la_{ri}^{(k)}$ by
$\la_{ri}^{(k)}+1$, and
the vector $\ze^{}_{\La}$ is considered to be zero,
if $\La(u)$ is not a pattern.
\ele

\bpf
The argument is based on Lemma~\ref{lem:relbc}.
As in the proof of Lemma~\ref{lem:arze},
separating the first
factor, we need to calculate the image
$B_r(-l^{\tss(k)}_{ri})\ts  C_s(v)\ts\eta$
for the respective value of $v$.
By \eqref{abc}, the operator
$B_r(u)$ commutes with $C_s(v)$ for $s\ne r$.
If $s=r$ then we consider two cases. If
$-l^{\tss(k)}_{ri}-v\ne 1$, then applying the first relation
of Lemma~\ref{lem:relbc} together with Lemma~\ref{lem:arze},
we find that
\ben
B_r(-l^{\tss(k)}_{ri})\ts  C_r(v)\ts\eta
=C_r(v)\ts B_r(-l^{\tss(k)}_{ri})\ts\eta
\een
and proceed by induction. If $v=-l^{\tss(k)}_{ri}-1$, then
we apply the second relation of
Lemma~\ref{lem:relbc} together with Lemma~\ref{lem:arze}
to get
\ben
B_r(-l^{\tss(k)}_{ri})\ts  C_r(-l^{\tss(k)}_{ri}-1)\ts\eta
=-A_{r+1}(-l^{\tss(k)}_{ri})
\ts A_{r-1}(-l^{\tss(k)}_{ri}-1)\ts\eta.
\een
One more application of Lemma~\ref{lem:arze} leads
to the desired formula.
\epf

The following theorem provides a basis of
the Gelfand--Tsetlin type
for the representation
$L(\la(u))$.

\bth\label{thm:basisgtylp}
The vectors $\ze_{\La}$ parameterized by all patterns $\La(u)$
associated with the highest weight $\la(u)$,
form a basis of the representation $L(\la(u))$
of the algebra $W(\pi)$.
\eth

\bpf
It is easy to verify that if the array of monic
polynomials obtained from $\La(u)$ by increasing
the entry $\la^{\tss(k)}_{ri}$ by $1$ is a pattern,
then the coefficient of the vector $\ze_{\La+\de_{ri}^{\tss(k)}}$
in the formula of Lemma~\ref{lem:brze} is nonzero.
This implies that each vector $\ze_{\La}\in L(\la(u))$
associated with a pattern $\La(u)$ is nonzero.

Furthermore, by Lemma~\ref{lem:arze},
$\ze_{\La}$ is an eigenvector
for all operators $A_r(u)$
with distinct sets of eigenvalues.
This shows that the vectors $\ze^{}_{\La}$
are linearly independent.

Finally, for each $i\in\{1,\dots,n\}$ and
$p_{i-1}< k\leqslant p_i$ the set of parameters
$\big(\la_{rj}^{(k)}\big)$
with $i\leqslant j\leqslant r\leqslant n$
forms a Gelfand--Tsetlin pattern associated with the highest weight
$\la^{(k)}$ of the irreducible representation
$L(\la^{(k)})$ of the Lie algebra $\gl_{q_k}$.
Hence, the
number of patterns $\La(u)$ coincides with the product
of dimensions $\dim L(\la^{(k)})$ for $k=1,\dots,l$.
Comparing this with \eqref{comdim}, we conclude that the number of
patterns coincides with $\dim L(\la(u))$.
\epf

Note that by Theorem~\ref{thm:basisgtylp},
we have the equality in \eqref{comdim}, and thus we recover
a result from \cite{bk:rs}
that the representation \eqref{tenprp} of $W(\pi)$
is irreducible.

\section{Action of the generators}

We will calculate the action of the generators of $W(\pi)$
in a normalized basis of $L(\la(u))$. For any pattern
$\La(u)$ associated to $\la(u)$ set
\begin{multline}
\non
N_{\La}=\prod_{(r,i)}\prod_{j=1}^{i-1}
\prod_{m=1}^{p_j}\prod_{k=1}^{p_i}
(l_j^{\tss(m)}-l_i^{\tss(k)})(l_j^{\tss(m)}-l_i^{\tss(k)}+1)
\dots (l_j^{\tss(m)}-l_{ri}^{\tss(k)}-1)\\
{}\times\prod_{j=i}^{r-1}\prod_{m=1}^{p_j}\prod_{k=1}^{p_i}
(l_{r-1,j}^{\tss(m)}-l_i^{\tss(k)})
(l_{r-1,j}^{\tss(m)}-l_i^{\tss(k)}+1)
\dots (l_{r-1,j}^{\tss(m)}-l_{ri}^{\tss(k)}-1),
\end{multline}
where the pairs $(r,i)$ run over the set of indices
satisfying the inequalities
$1\leqslant i\leqslant r\leqslant n-1$.
This constant is clearly nonzero for any pattern $\La(u)$.
Introduce normalized vectors $\xi_{\La}\in L(\la(u))$ by
\ben
\xi_{\La}=N_{\La}^{-1}\ts\ze_{\La}.
\een
By Theorem~\ref{thm:basisgtylp},
the vectors $\xi_{\La}$
form a basis of the representation $L(\la(u))$.
The algebra $W(\pi)$ is generated by the coefficients
of the polynomials $A_r(u)$ with $r=1,\dots,n$ and
the coefficients of
the polynomials $B_r(u)$ and $C_r(u)$ with $r=1,\dots,n-1$.
Since $B_r(u)$ and $C_r(u)$
are polynomials in $u$ of degree
less than $p_1+\dots+p_r$, it suffices to
find the values of these polynomials at $p_1+\dots+p_r$
different values
of $u$. The polynomial can then be calculated by the Lagrange
interpolation formula. For these values we take the numbers
$-l^{\tss(k)}_{ri}$ with $i=1,\dots,r$ and $k=1,\dots,p_i$.

\bth\label{thm:dgactylp}
We have
\beql{aractba}
A_r(u)\ts \xi^{}_{\La}=\la_{r1}(u)\dots\la_{rr}(u-r+1)\ts \xi^{}_{\La},
\eeq
for $r=1,\dots,n$, and
\begin{align}\label{bractba}
B_r(-l^{(k)}_{ri})\ts \xi^{}_{\La}&=
-\la_{r+1,1}(-l^{\tss(k)}_{ri})\dots\la_{r+1,r+1}(-l^{\tss(k)}_{ri}-r)
\ts \xi^{}_{\La+\de_{ri}^{(k)}},\\
C_r(-l^{(k)}_{ri})\ts \xi^{}_{\La}&=
\la_{r-1,1}(-l^{\tss(k)}_{ri})\dots\la_{r-1,r-1}(-l^{\tss(k)}_{ri}-r+2)
\ts \xi^{}_{\La-\de_{ri}^{(k)}},
\non
\end{align}
for $r=1,\dots,n-1$, where $\xi^{}_{\La\pm\de_{ri}^{(k)}}$
corresponds to the pattern
obtained from $\La(u)$ by replacing $\la_{ri}^{(k)}$ by
$\la_{ri}^{(k)}\pm 1$.
\eth

\bpf
The formulas for the action of $A_r(u)$ and $B_r(-l^{(k)}_{ri})$
follow respectively from Lemmas~\ref{lem:arze} and \ref{lem:brze}
by taking into account the normalization constant.
Now consider the vector $C_r(-l^{(k)}_{ri})\ts\xi^{}_{\La}$.
Arguing as in the proof of Lemma~\ref{lem:arze},
and using \eqref{aractba},
we find that
\ben
\bal
A_s(u)\ts C_r(-l^{(k)}_{ri})\ts\xi^{}_{\La}
{}&=C_r(-l^{(k)}_{ri})\ts A_s(u)\ts\xi^{}_{\La}\\
{}&=\la_{s1}(u)\dots\la_{ss}(u-s+1)
\ts C_r(-l^{(k)}_{ri})\ts\xi^{}_{\La}
\eal
\een
for $s\ne r$, while
\ben
\bal
A_r(u)\ts C_r(-l^{(k)}_{ri})\ts\xi^{}_{\La}
{}&=\frac{u+l^{(k)}_{ri}-1}{u+l^{(k)}_{ri}}
C_r(-l^{(k)}_{ri})\ts A_r(u)\ts\xi^{}_{\La}\\
{}&=\frac{u+l^{(k)}_{ri}-1}{u+l^{(k)}_{ri}}
\ts\la_{r1}(u)\dots\la_{rr}(u-r+1)
\ts C_r(-l^{(k)}_{ri})\ts\xi^{}_{\La}.
\eal
\een
If $\la_{ri}^{(k)}=\la_{r+1,i+1}^{(k)}$, then
the vector $\xi^{}_{\La-\de_{ri}^{(k)}}$ is zero and we need
to show that
$C_r(-l^{(k)}_{ri})\ts \xi^{}_{\La}=0$.
Indeed, otherwise the vector
$C_r(-l^{(k)}_{ri})\ts \xi^{}_{\La}$
must be proportional
to a certain basis vector of $L(\la(u))$.
However, this is impossible because none
of the basis vectors has the same set of eigenvalues
as $C_r(-l^{(k)}_{ri})\ts \xi^{}_{\La}$.

If $\la_{ri}^{(k)}-\la_{r+1,i+1}^{(k)}\geqslant 1$,
then by the same argument we have
\ben
C_r(-l^{(k)}_{ri})\ts \xi^{}_{\La}
=\al\ts \xi^{}_{\La-\de_{ri}^{(k)}}
\een
for a certain constant $\al$. Its value
is found by the application of the operator
$B_r(-l^{(k)}_{ri}+1)$ to the vectors on both sides
with the use of \eqref{aractba},
\eqref{bractba} and the second relation in
Lemma~\ref{lem:relbc}.
\epf

Note that in the particular case of a rectangular pyramid $\pi$
the normalized basis $\{\xi_{\La}\}$ coincides with
the basis of \cite{m:gt} constructed in a different way.

Let us denote by $\pi'$ the pyramid with the rows $p_1,\dots,p_{n-1}$.
Then the finite $W$-algebra
$W(\pi')$ may be identified with the subalgebra
of $W(\pi)$ generated by the elements \eqref{gener},
excluding all $h_{n}^{(r)}$, $e_{n-1}^{(r)}$ and $f_{n-1}^{(r)}$.
Theorem~\ref{thm:dgactylp} implies the
following {\it branching rule\/} for the reduction
$W(\pi)\downarrow W(\pi')$ and thus shows that the basis
$\{\xi_{\La}\}$ is consistent with the chain of subalgebras
\eqref{chainw}.

\bco\label{cor:braylp}
The restriction of the $W(\pi)$-module $L(\la(u))$
to the subalgebra $W(\pi')$ is isomorphic to
the direct sum
of irreducible highest weight $W(\pi')$-modules $L'(\mu(u))$,
\ben
L(\la(u))|^{}_{W(\pi')}
\cong\underset{\mu(u)}{\scloplus}\ts L'(\mu(u)),
\een
where $\mu(u)$ runs over all $(n-1)$-tuples of monic
polynomials in $u$ of the form
$\mu(u)=\big(\mu_1(u),\dots,\mu_{n-1}(u)\big)$,
such that
\ben
\mu_i(u)=(u+\mu^{(1)}_i)\ts(u+\mu^{(2)}_i)
\dots (u+\mu^{(p_i)}_i),\qquad i=1,\dots,n-1,
\een
and the following conditions hold:
\ben
\la_{i}^{(k)}-\mu_i^{(k)}\in\ZZ_+\qquad\text{and}\qquad
\mu_i^{(k)}-\la_{i+1}^{(k)}\in\ZZ_+
\een
for $k=1,\dots,p_i$ and
$1\leqslant i\leqslant r\leqslant n-1$.
\qed
\eco

For each $i=1,\dots,n-1$ introduce the polynomials
$\tau_{ni}(u)$ and $\tau_{in}(u)$
with coefficients in $W(\pi)$ by
the formulas
\beql{lowraise}
\bal
\tau_{ni}(u)&=C_{n-1}(u)\ts C_{n-2}(u)\dots C_i(u),\\
\tau_{in}(u)&=B_i(u)\ts B_{i+1}(u)\dots B_{n-1}(u).
\eal
\eeq
Define the vector
$\ze_{\mu}\in L(\la(u))$
corresponding to the $(n-1)$-tuple of polynomials
$\mu(u)=(\mu_1(u),\dots,\mu_{n-1}(u))$
by the formula
\ben
\ze_{\mu}=\prod_{i=1}^{n-1}\prod_{k=1}^{p_i}
\Big(\tau_{ni}(-m_i^{(k)}-1)\dots
\tau_{ni}(-l_i^{(k)}+1)
\ts\tau_{ni}(-l_i^{(k)})\Big)\ts\ze,
\een
where the ordering of the factors corresponds to
increasing indices $i$ and $k$, and we used the notation
\ben
m_i^{(k)}=\mu_i^{(k)}-i+1\qquad\text{and}\qquad
l_i^{(k)}=\la_i^{(k)}-i+1.
\een
Due to Theorem~\ref{thm:dgactylp}, the vector
$\ze_{\mu}$ generates a $W(\pi')$-submodule
of $L(\la(u))$, isomorphic to $L'(\mu(u))$.
Moreover, the operators $\tau_{ni}(-m_i^{(k)})$
and $\tau_{in}(-m_i^{(k)})$ take $\ze_{\mu}$
to the vectors proportional to
$\ze_{\mu-\de_i^{(k)}}$ and $\ze_{\mu+\de_i^{(k)}}$,
respectively.
So, the polynomials \eqref{lowraise} valued
at appropriate points can be regarded
as the {\it lowering\/} and {\it raising\/} operators
for the reduction $W(\pi)\downarrow W(\pi')$;
cf. \cite[Chapter~5]{m:yc}.

\section*{Acknowledgments}

The authors acknowledge the support of
the Australian Research Council.
The first author is supported in part
by the CNPq grant (processo 307812/2004-9)
and by the Fapesp grant
(processo 2005/60337-2).
He is grateful to
the University of Sydney for the warm hospitality
during his visit.

\end{document}